\documentclass[a4paper, 11pt]{article}

\usepackage[
            includefoot,  
            marginparwidth=0in,     
            marginparsep=0in,       
            margin=1in,               
            includemp]{geometry}

\usepackage{bbm}
\usepackage{float}
\usepackage{graphicx}                  
\usepackage{amssymb}
\usepackage{amsfonts}
\RequirePackage{amsmath}
\RequirePackage{amsthm}

\usepackage{color}

\usepackage{fancyhdr}

\newtheorem{thm}{Theorem}[section]
\newtheorem{defn}[thm]{Definition}
\newtheorem{prop}[thm]{Proposition}
\newtheorem{lem}[thm]{Lemma}

\newtheorem{cor}[thm]{Corollary}

\DeclareMathOperator{\Var}{Var}
\DeclareMathOperator{\cov}{cov}

\newcommand{\R}{\mathbb{R}}

\newcommand{\bS}{\mathbb{S}}

\begin{document}

\title{Growth estimate on optimal transport maps via concentration inequalities}
\author{Max Fathi\footnote{Université Paris Cité and Sorbonne Université, CNRS, Laboratoire Jacques-Louis Lions and Laboratoire de Probabilit\'es, Statistique et Mod\'elisation, F-75013 Paris, France\\
and DMA, École normale supérieure, Université PSL, CNRS, 75005 Paris, France \\
and Institut Universitaire de France\\ mfathi@lpsm.paris }}
\date{\today}

\maketitle
{\centering\footnotesize \textit{Dedicated to Patrick Cattiaux, with gratitude for the advice over the years.}\par}
\abstract{We give an alternative proof and some extensions of results of \cite{CFS} on polynomial upper bounds on the Brenier map between probability measures under various conditions on the densities. The proofs are based on the monotonicity of the map and various concentration inequalities, as already used in \cite{CF} to prove quadratic growth for transport maps from the standard Gaussian onto log-concave measures.
 
\section{Introduction}

In this work, we are interested in the asymptotic growth of the $L^2$ optimal transport (or Brenier map) between probability measures under growth conditions on the densities. This line of research is motivated by the classical Caffarelli contraction theorem \cite{Caf}, which states that the optimal transport map from a standard Gaussian measure onto a uniformly log-concave measure is globally lipschitz, with an explicit estimate that is dimension-free. There have been various alternative proofs and generalizations of this theorem, see for example \cite{Kol, ChPo, FGP, KM}. 

In situations where we cannot prove an analogue of Caffarelli's theorem, it is natural to seek to prove weaker statements. One example is to study whether the map satisfies a linear growth bound of the form 
$$|T(x)| \leq C(1 + |x|),$$
ideally with an explicit constant prefactor. This is obviously weaker than proving that $T$ is globally Lipschitz. This question was first studied in \cite{CF} for maps from a standard Gaussian onto a log-concave measure and also pursued in a much more general setting in \cite{CFS} (whose results we shall detail later). In both \cite{CF, CFS}, the growth estimates where then bootstrapped into local Lipschitz estimates for the optimal map under growth assumptions on the second derivatives of the log-densities of the measures, but here we shall not pursue that second line of inquiry, and only focus on growth estimates. 

The methods of \cite{CF} and \cite{CFS} for proving growth estimates are quite different: \cite{CF} relies on the monotonicity of the optimal map and uses a concentration inequality for the target measure to deduce a growth estimate, while \cite{CFS} relies on the Monge-Amp\`ere PDE satisfies by the map, and a maximum principle-type argument. The purpose of this work is to show that the method of \cite{CF} can be generalized to prove the growth estimates of \cite{CFS}, as well as several extensions. 

In all the sequel, $C$ will be a positive constant that may change from line to line, and depends on parameters whose influence we do not precisely quantity in the desired statement. It shall always be computable. 

\section{The general setup}

Let us first recall the Brenier-McCann \cite{Bre, McC} theorem on existence and uniqueness of monotone transport maps: 

\begin{thm}
Let $\mu \in \mathcal{P}(\R^d)$ be absolutely continuous with respect to the Lebesgue measure, and let $\nu \in \mathcal{P}(\R^d)$. There exists a convex function $\varphi$ such that $T = \nabla \varphi$ is a transport map from $\mu$ onto $\nu$. Moreover, $T$ is the unique (up to null sets) transport map of this form. 
\end{thm}

The density assumption on $\mu$ can be slightly relaxed, but we shall only consider this situation here. We shall be interested here in growth estimates at infinity on $T(x)$ under various conditions on $\mu$ and $\nu$.  

The analysis of \cite{CFS} relies on the fact that $\varphi$ is the unique solution (up to constants) to the Monge-Amp\`ere equation 
$$\mu(x) = \nu(\nabla \varphi(x))\det(\nabla^2 \varphi(x))$$
where we identify the two measures and their density. We refer to \cite{Fig} for the qualitative theory of regularity for solutions to this PDE. The proofs in \cite{CFS} are based on the maximum principle applied to functionals of the solution, as pioneered in \cite{Caf}. 

We shall instead use the monotonicity of the optimal transport map to restrict the regions where points are sent compared to each other, and then use concentration inequalities to bound the mass of those regions in a way that depends on $|T(x)|$. The following is an argument of \cite{CF}

\begin{lem} \label{key_lemma}
Let $T = \nabla \varphi$ with $\varphi$ convex. Then for any $u \in \bS^{d-1}$ and $\lambda > 0$ we have
$$T(B(x + 2\lambda u, \lambda) \subset T(x) + \{v; \hspace{3mm} \langle v, u \rangle \geq -|v|/2\}.$$
In other words, setting 
$$f(z) = \frac{2}{3}\langle z - T(x), u \rangle + \frac{|z-T(x)|}{3},$$ 
we have
$$T(B(x + 2\lambda u, \lambda) \subset \{f \geq 0\}.$$
\end{lem}

\begin{proof}
Since $\varphi$ is convex, we have for any $y$
$$\langle T(y) - T(x), y - x \rangle \geq 0.$$
Letting $v = T(y) - T(x)$ with $y \in B(x + 2\lambda u, \lambda)$, there exists $w \in B(0,1)$ such that
$$\langle v, 2\lambda u + \lambda w \rangle \geq 0.$$
This implies
$$\langle v, u \rangle \geq - |v|/2$$
which concludes the proof. 
\end{proof}

The following statement converts the previous lemma into a bound on $|T(x)|$. 

\begin{thm} \label{main_abstract_thm}
Let $\psi$ be a continuous decreasing function such that $\psi(r) \geq \nu(|z| \geq r)$ for all $r \geq r_0$. Then for any $\lambda > 0$ and $x \in \R^d$, we have
$$|T(x)| \leq \max(3r_0,  3\psi^{-1}(\mu(B(x + 2\lambda u, \lambda)))$$
with $u = T(x)/|T(x)|$. 
\end{thm}

\begin{proof}
We have
\begin{align*}
&\left\{\frac{2}{3}\langle z, u\rangle - \frac{2}{3}|T(x)| + \frac{1}{3}|z-T(x)|\geq 0\right\} \\
&\subset \left\{\frac{2}{3}\langle z, u\rangle - \frac{2}{9}|T(x)| \geq 0\right\} \cup \left\{- \frac{4}{9}|T(x)| + \frac{1}{3}|z-T(x)|\geq 0\right\} \\
&\subset \{|z| \geq |T(x)|/3\} \cup \{|z| \geq |T(x)|/3\} \\
&= \{|z| \geq |T(x)|/3\}.
\end{align*}
Therefore, if $|T(x)| \geq 3r_0$, applying Lemma \ref{key_lemma} yields
$$\nu(T(B(x + 2\lambda u, \lambda)) \leq \psi(|T(x)|/3).$$
However, since $T$ is a transport map, 
$$\nu(T(B(x + 2\lambda u, \lambda)) = \mu(B(x + 2\lambda u, \lambda)).$$
Therefore, 
$$\mu(B(x + 2\lambda u, \lambda)) \leq \psi(|T(x)|/3).$$
\end{proof}

Using Theorem \ref{main_abstract_thm} requires being able to bound the probability (under the target measure) of taking large values, i.e. a concentration inequality. To connect with the literature, we think it is convenient to go through general notions of concentration inequalities, defined as following: 

\begin{defn}
Let $r _0 \geq 0$ and let $\varphi$ be a decreasing function from $\R_+$ to $[0,1]$. A probability measure $\nu$ satisfies a concentration inequality with parameters $(\varphi, r_0)$ if for any $1$-Lipschitz function $f$ and any $r \geq r_0$ we have
$$\nu\left(\left\{f \geq \int{fd\nu} + r \right\}\right) \leq \varphi(r).$$
\end{defn}

In principle the function $\varphi$ should be merely non-increasing, but in practice we will only use functions that are decreasing, so we take it as part of the assumptions to simplify some notations. 

\begin{cor} \label{main_abstract_cor}
Assume $\nu$ is centered, satisfies a concentration inequality with parameters $(\varphi, r_0)$, and that $\int{|z|d\nu} \leq M$. Then for any $\lambda > 0$ and $x \in \R^d$, we have
$$|T(x)| \leq \max(M + 3r_0,  M + 3\varphi^{-1}(\mu(B(x + 2\lambda u, \lambda))).$$
\end{cor}

This can be viewed as a Corollary of Theorem \ref{main_abstract_thm}, by using the fact that the norm is $1$-Lipschitz, but we go back to Lemma \ref{key_lemma} to give a proof, in order to highlight the structure (and slightly improve a constant by the way). 

\begin{proof}
Let $x \in \R^d$, and assume without loss of generality that $|T(x)| \geq M + 3r_0$

As in the proof of Lemma \ref{key_lemma}, let
$$f(z) = \frac{2}{3}\langle z - T(x), u \rangle + \frac{|z-T(x)|}{3}$$ 
and take $u = T(x)/|T(x)|$. Then $f$ is a $1$-lipschitz function, so for any $r \geq r_0$, 
$$\nu\left(\left\{f \geq \int{fd\nu} + r\right\}\right) \leq \varphi(r).$$
We have
\begin{align*}
\int{fd\nu} &= -\frac{2}{3}|T(x)| + \frac{1}{3}\int{|z-T(x)|d\nu(z)} \\
&\leq - |T(x)|/3 + M/3.
\end{align*}
Hence taking $r = (|T(x)| - M)/3$, which is greater than $r_0$, and applying the concentration inequality yields
$$\nu\left(\left\{f \geq 0\right\}\right) \leq \varphi((|T(x)| - M)/3).$$
Therefore, applying Lemma \ref{key_lemma} yields
$$\nu(T(B(x + 2\lambda u, \lambda)) \leq \varphi((|T(x)| - M)/3).$$
However, since $T$ is a transport map, 
$$\nu(T(B(x + 2\lambda u, \lambda)) = \mu(B(x + 2\lambda u, \lambda)).$$
Therefore, 
$$\mu(B(x + 2\lambda u, \lambda)) \leq \varphi((|T(x)| - M)/3).$$
Recalling that $u = T(x)/|T(x)|$, this concludes the proof. 
\end{proof}

\section{Concentration estimates}

\subsection{Concentration inequalities}

In this section we present some typical families of concentration inequalities, and various (mostly classical) sufficient conditions for ensuring they hold. 

\begin{defn}
A probability measure $\mu$ is $\sigma^2$-subgaussian if for any $1$-lipschitz function $f$ with $\int{fd\mu} = 0$ we have
$$\int{\exp(\lambda f)d\mu} \leq \exp(\sigma^2\lambda^2/2).$$
Applying Markov's inequality and optimizing in $\lambda$, in particular subgaussian measures satisfy
$$\mu\left(\left\{f \geq \int{fd\mu} + r\right\}\right) \leq \exp\left(-\frac{r^2}{2\sigma^2}\right)$$
for any $1$-lipschitz function $f$ and any $r > 0$. 
\end{defn}
The standard Gaussian measure is $1$-subgaussian. Other examples include measures satisfying a logarithmic Sobolev inequality or a Talagrand transport-entropy inequality, log-bounded or log-lipschitz perturbations of uniformly convex measures, as well as compactly supported measures. 

The next proposition summarizes classical results on measures satisfying the Gaussian concentration property: 

\begin{prop}
Let $W$ be such that $\nu = e^{-W}$ is a probability density on $\R^d$. 

\begin{enumerate}
\item If $\nabla^2 W \geq \kappa $ with $\kappa > 0$, then $\nu$ is $\kappa$-subgaussian;

\item If $W = W_1 + W_2$ with $\nabla^2 W_1 \geq \kappa $ and $||W_2||_\infty \leq \delta$ then $\nu$ is $\kappa \exp(-2\delta)$-subgaussian;

\item If $W = W_1 + W_2$ with $\nabla^2 W_1 \geq \kappa $ and $||W_2||_{lip} \leq \delta$ then $\nu$ is $\sigma^2$-subgaussian, for some $\sigma^2$ depending only on $\kappa$ and $\delta$ (but not on $d$);
\end{enumerate}
\end{prop}

\begin{proof}
In all three situations, a logarithmic Sobolev inequality (LSI) holds, which then implies the Gaussian concentration property by the Herbst argument \cite[Proposition 5.4.1]{BGL}. 

For part (i), the LSI is the outcome of the classical Bakry-Emery theorem \cite{BE85}. For (ii), it is a combination of the LSI in case (i) and the Holley-Stroock perturbation lemma \cite[Proposition 5.1.6]{BGL}. Finally, (iii) holds as a consequence of the Aida-Shigekawa theorem on stability of the LSI for log-Lipschitz perturbations \cite{AS}. 
\end{proof}

Another important type of concentration inequality is exponential concentration: 

\begin{defn}
A probability measure satisfies an exponential concentration inequality with parameters $(c,\sigma)$ if for any $1$-Lipschitz function we have
$$\mu\left(f \geq \int{fd\mu} + r\right) \leq c\exp(-r/\sigma) \hspace{3mm} \forall r \geq 0.$$
\end{defn}

An important class of measures satisfying an exponential concentration inequality are measures satisfying a Poincar\'e inequality. In particular, it covers the situation  of log-concave measures, as proved in \cite{Kla23}:  

\begin{prop}
Let $\nu = e^{-W}$ be a probability measure on $\R^d$, and assume that $W$ is convex and that $||\cov(\mu)||_{op} \leq 1$. Then $\mu$ satisfies an exponential concentration inequality with parameters $(c_1, c_2\sqrt{\log d})$, with $c_1, c_2$ universal constants.  
\end{prop}

It is conjectured that the factor $\sqrt{\log d}$ is unnecessary. We refer to the lecture notes \cite{KL} for background on this conjecture.  

As for polynomial concentration bounds, we have: 

\begin{prop} Let $W$ be a positive function such that $\nu = W^{-d}$ is a probability measure on $\R^d$. If there exists $p > 1$ and $\beta > 1$ such that
\begin{equation}
\liminf_{|y|\rightarrow + \infty} \frac{W(y)}{|y|^p} > 0; \hspace{3mm} \liminf_{|y|\rightarrow + \infty} \frac{y \cdot \nabla W(y)}{W(y)} \geq \beta
\end{equation}
then for any $\ell < d(p-1)$, $\nu$ satisfies a concentration inequality with parameters $(r_0, Cr^{-\ell})$, for some computable constants $r_0$ and $C$. 
\end{prop}

We shall not actually use this Proposition, since it would give a slightly sub-optimal exponent, but we state  it to highlight that the assumptions used in \cite[Theorem 2.2]{CFS} are indeed the same as those used in known results for proving polynomial concentration inequalities. One could also get similar concentration inequalities for measures with density of the form $W(x)^{-d-\alpha}$ with $\alpha > 0$ and $W$ convex, see \cite[Remark 4.12]{CGGR}. 

\begin{proof}
We shall use the Lyapunov function method for proving concentration inequalities, as proposed in \cite{CGGR}. Let $L = \Delta - d \nabla \log W \cdot \nabla$ be the generator of the reversible drifted Brownian motion with invariant probability measure $\nu$. Setting $g(x) = 1 + |x|^k/k$ with $k > 2$, we have
$$Lg(x) = d(k-1)|x|^{k-2} - d|x|^{k-2}\frac{x \cdot \nabla W(x)}{W(x)}.$$
If we take $R \geq 2$ and such that $\frac{y \cdot \nabla W(y)}{W(y)}\geq \beta'$ on $B(0,R)^c$, with $\beta' < \beta$, and $\alpha = -\inf_{B(0,R)} \frac{y \cdot \nabla W(y)}{W(y)}$, then using that for $|x| \geq 2$ we have $|x|^k \geq g(x)$, we have
\begin{align*}
Lg(x) &\leq d(k-1 + \alpha)R^{k-2} \mathbbm{1}_{B(0,R)} -d(\beta' - (k-1))|x|^{k-2}\mathbbm{1}_{B(0,R)^c} \\
&\leq  d(k-1 + \alpha)R^{k-2} \mathbbm{1}_{B(0,R)} -d(\beta' - (k-1))g(x)^{(k-2)/k}\mathbbm{1}_{B(0,R)^c} \\
&\leq d(2(k-1) + \alpha + \beta')R^{k-2} \mathbbm{1}_{B(0,R)} -d(\beta' - (k-1))g(x)^{(k-2)/k}.
\end{align*}
Hence as long as $\beta > k-1$, we can find $R, C_1, C_2 > 0$ such that
\begin{equation}Lg(x) \leq C_1 \mathbbm{1}_{B(0,R)} - C_2g(x)^{(k-2)/k} \hspace{3mm} \forall x \in \R^d.
\end{equation}
According to \cite[Theorem 2.8]{CGGR}, this implies that there is a (computable) constant $C > 0$ such that the weighted Poincar\'e inequality
\begin{equation}
\Var_{\mu}(f) \leq C\int{(1+|x|^2)|\nabla f|^2d\mu}.
\end{equation}
The weight arises as $1 + |\nabla g|^2/g^{2(k-2)/k}$. Applying \cite[Theorem 2.8]{CGGR} requires a local Poincar\'e inequality for $\nu$, but this assumption is easily satisfied for measures with a positive continuous density (inducing a constant that depends on $L^\infty$ bounds on the density within a large ball). Note that the value of $k > 2$ does not appear explicitly anymore. Its choice may have affected the value of the constant, but tuning it does not seem to help identify better assumptions. 

According to \cite[Corollary 4.2]{BoLe07}, this weighted Poincar\'e inequality implies a concentration inequality with parameters $t_0 > 0$ and $\phi(r) = Cr^{-\ell}$, with computable constants $t_0$ and $C$, for any $\ell$ such that
$$\int{(1+|x|^2)^{\ell/2}d\mu} < \infty.$$
Using that $W(x) \geq C(1 + |x|^p)$, we see that $\int{(1+|x|^2)^{\ell/2}d\mu} < \infty$ as soon as $\ell < d(p-1)$. 
\end{proof}

\subsection{Ball probability lower bounds}

\begin{lem}[Ball probability lower bounds for polynomial-like densities] \label{lem_ball_probab_polynom}
Let $V$ be a positive function on $\R^d$ such that $V^{-d}$ is a probability density. 
\begin{enumerate}
\item Assume that $V(z) \leq L(1 + |z|^q)$ for some $q > 1$ and $L > 0$. Then for any $u \in \mathbb{S}^{d-1}$ we have
$$\mu(B(x + 4|x|u, 2|x|)) \geq \alpha(1 + |x|)^{-(q-1)d} \hspace{3mm} \forall x \in \R^d$$
for a constant $\alpha$ that depends on $d, q, L$ and $V$, but not on $x$ nor $u$. 

\item Assume that  
$$|\nabla \log V| \leq A(1 + |x|)^{-1}.$$
Then $$\mu(B(x,1/2)) \geq \exp(-3Ad/2)(1+2|x|)^{-2Ad}\mu(B(0,1/2)) \hspace{3mm} \forall x \in \R^d.$$
\end{enumerate}
\end{lem}

\begin{proof}
Let us start with proving 1. Without loss of generality assume $|x| \geq 1$. By assumption, we have
\begin{align*}
\mu(B(x + 4|x|u, 2|x|)) &\geq L^{-d}\int_{B(x + 4|x|u, 2|x|)}{(1+|z|^q)^{-d}dz} \\
&\geq C\int_{B(x + 4|x|u, 2|x|)}{|z|^{-qd}dz} \\
&\geq C|x|^{d-qd}\int_{B(x/|x| + 4u, 2)}{|y|^{-qd}dy}\\
&\geq C(1 + |x|)^{-(q-1)d},
\end{align*}
where we used the change of variable $z = |x|y$. Note that for any $x \neq 0$ and $u$ we have 
$$B(x/|x| + 4u, 2) \subset  B(0,7)\backslash B(0,1),$$
so that the integral on the fourth line is over a subset of a region uniformly bounded away from both the origin and infinity. The constant $C$ depends on bounds on $V$ in a bounded region, $L$, $p$ and $d$, but not on $x$ nor $u$.

We shall now prove 2. We will first prove a ratio bound for $V$, and then conclude by a change of variable. 

Let $y \in B(0,1/2)$ and $|x| > 1$ then, using the assumption on $V$ we have
\begin{align*}
\frac{V(x)}{V(y)} &= \exp\left(\int_0^1{\nabla \log V(tx + (1-t)y) \cdot (x-y) dt}\right) \\
&\leq \exp\left(A|x-y|\int_0^1{(1 + |tx + (1-t)y|)^{-1}dt}\right) \\
&\leq \exp\left(2A|x|\int_0^1{(1 + t|x| + (1-t)|y|)^{-1}dt}\right) \\
&\leq \exp\left(2A|x|\int_0^1{(1/2 + t|x|)^{-1}dt}\right) \\
&\leq \exp\left(4A|x|\int_0^1{(1 + 2t|x|)^{-1}dt}\right) \\
&\leq (1+2|x|)^{2A}.
\end{align*}

If $|x| \leq 1$, since $\log V$ is $A$-Lipschitz, we have
\begin{equation*}
\frac{V(x)}{V(y)} \leq \exp\left(A|x-y|\right) \leq \exp(3A/2)
\end{equation*}
So for any $x \in \R^d$ and $y \in B(0,1/2)$ we have
$$\frac{V(x)}{V(y)} \leq \exp(3A/2)(1+2|x|)^{2A}.$$
Therefore,
\begin{align*}
\mu(B(x,1/2)) &= \int_{B(x,1/2)}{V(z)^{-d}dz} \\
&= \int_{B(0,1/2)}{V(x+z)^{-d}dz} \\
&\geq \int_{B(0,\sqrt{d})}{V(z)^{-d}\exp(-3Ad/2)(1+2|x|)^{-2Ad}dz} \\
&= \exp(-3Ad/2)(1+2|x|)^{-2Ad}\mu(B(0,1/2)).
\end{align*}
\end{proof}

\section{Growth estimates on optimal transport maps}

The following is a generalization of \cite[Theorem 2.4]{CFS}: 

\begin{thm} \label{thm_growth_polyn_to_subgaussian}
Assume that the starting measure $\mu$ has density $V^{-d}$ with $V$ satisfying 
$$|\nabla \log V| \leq A(1 + |x|)^{-1}$$
and that the target measure $\nu$ is centered and $\sigma^2$-subgaussian. 
Then the Brenier map satisfies 
$$|T(x)| \leq \sqrt{d}\sigma\left(3 + 3\sqrt{2A \log 2 + 5A - d^{-1}\log \omega_d + \log V(0)  + 2A\log(1+ |x|)}\right)$$
where $\omega_d$ is the volume of the unit ball in dimension $d$. 
\end{thm}

Note that in high dimension, $- d^{-1}\log \omega_d$ is equivalent to $d\log(d)/2$. 

Remark that under these assumptions, $\nu$ might not have a density, and the transport map might not be continuous. In particular, we cannot directly use the Monge-Ampère PDE in a pointwise sense to analyze this situation. 

\begin{proof}[Proof of Theorem \ref{thm_growth_polyn_to_subgaussian}]
We shall apply Corollary \ref{main_abstract_cor} with $\lambda = 1/2$. The Gaussian concentration assumption allows to take $r_0 = 0$ and $\varphi(r) = \exp(-r^2/(2\sigma^2))$. Moreover, since the measure is centered, we have
\begin{align*}
\frac{1}{2}\int{x_i^2d\nu} &= \lim_{t \rightarrow 0}t^{-2}\int{\exp(tx_i) - 1 - tx_1d\nu}
\leq \lim_{t \rightarrow 0}t^{-2}(\exp(t^2\sigma^2/2) - 1) = \sigma^2/2.
\end{align*}
Summing over all coordinates and using Jensen's inequality, we see we can take $M =  \sqrt{d}\sigma$. 

Using part 2 of Lemma \ref{lem_ball_probab_polynom}, denoting $u = T(x)/|T(x)|$, we have 
\begin{align*}
\mu(B(x + u,1/2)) &\geq \exp(-3Ad/2)(1+2|x+u|)^{-2Ad}\mu(B(0,1/2))\\
& \geq 2^{-2Ad}\exp(-3Ad/2)(1+|x|)^{-2Ad}\mu(B(0,1/2)).
\end{align*}
Applying Theorem \ref{main_abstract_thm}, we get 
$$|T(x)| \leq \sqrt{d}\sigma\left(1 + 3\sqrt{2A \log 2 + 3A - d^{-1}\log \mu(B(0,1/2)) + 2A\log(1+ |x|)}\right).$$
Finally, since $\log V$ is $A$-Lipschitz, we have
\begin{align*}
\mu(B(0,1/2)) &\geq \int_{B(0,1/2}{e^{-d\log V(0) - Adx}dx} 
&\geq \exp(-d\log V(0) -d(A + 1/2 + \log 2))\omega_d
\end{align*}
where $\omega_d$ is the volume of the unit ball in dimension $d$. Using this bound (and simplifying a bit the numerical constants) concludes the proof. 
\end{proof}

With the same arguments, we can get a similar estimates for measure satisfying an exponential concentration inequality: 

\begin{thm} \label{thm_growth_polyn_to_exponential}
Assume that the starting measure $\mu$ has density $V^{-d}$ with $V$ satisfying 
$$|\nabla \log V| \leq A(1 + |x|)^{-1}$$
and that the target measure $\nu$ satisfies an exponential concentration inequality with parameters $(c,\sigma)$. 
Then the Brenier map satisfies 
$$|T(x)| \leq 2\sqrt{cd}\sigma + 3\sigma d\left(d^{-1}\log c + \log V(0) + 2A \log 2 + 3A + 2 - d^{-1}\log \omega_d + 2A\log(1 + |x|) \right).$$
\end{thm}

The proof is essentially the same as for Theorem \ref{thm_growth_polyn_to_subgaussian}, so we omit it. The moment bound $\int{|x_i|d\nu} \leq 2\sqrt{c}\sigma$ is obtained simply by integrating the concentration inequality. 

For the particular case of log-concave measures, given the current known concentration inequalities, we get the following: 
\begin{cor}
Assume that the starting measure $\mu$ has density $V^{-d}$ with $V$ satisfying 
$$|\nabla \log V| \leq A(1 + |x|)^{-1}$$
and that the target measure $\nu$ is isotropic log-concave. 
Then the Brenier map satisfies 
$$|T(x) | \leq Cd(\log d)^{1/2}(1 + \log d + \log(1 + |x|))$$
where the constant $C$ only depends on $A$ and $V(0)$. 
\end{cor}

We now consider target measures with polynomially decaying densities. We get the following generalization of \cite[Theorem 2.2]{CFS}.

\begin{thm}
Consider the Brenier map from $\mu = V^{-d}dx$ onto $\nu = W^{-d}dx$, with both $V$ and $W$ continuous and positive. Assume moreover that there exists $p, q > 1$ and $L > 0$ such that
$$V(x) \leq L(1 + |x|^q); \hspace{3mm} W(y) \geq M(1 + |y|^p) \hspace{2mm} \forall x, y \in \R^d.$$
Then there is a computable constant $C > 0$ (that depends on $p,q,d,L$ and $M$ and the behaviour of $V$ in a fixed bounded region) such that
$$|T(x)| \leq C(1 + |x|)^{\frac{q-1}{p-1}} \hspace{3mm} \forall x \in \R^d.$$
\end{thm}

Compared to \cite[Theorem 2.2]{CFS}, we eliminate some assumptions on the exponents and the densities, answering in particular the question raised in \cite[Remark 2.3]{CFS}. 

\begin{proof}
We shall apply Theorem \ref{main_abstract_thm}. For the starting measure, applying Lemma \ref{lem_ball_probab_polynom} (part 1) yields
$$\mu(B(x + 4|x|u, 2|x|)) \geq \alpha(1 + |x|)^{-(q-1)d} \hspace{3mm} \forall x \in \R^d$$
with $u = T(x)/|T(x)|$. 

For the target measure, we have 
\begin{align*}
\nu(|y| \geq r) &= \int_{B(0, r)^c}{W(y)^{-d}dy} \\
&\leq \int_{B(0, r)^c}{M^{-d}(1 + |y|^p)^{-d}dy} \\
&\leq C\int_r^\infty{t^{-d(p-1)-1}dt} \\
&\leq Cr^{-d(p-1)}.
\end{align*}

Combining the two estimates via Theorem \ref{main_abstract_thm} concludes the proof. 
\end{proof}

\hspace{3mm}

\textbf{\underline{Acknowledgments}}: This work has received support under
the program "Investissement d'Avenir" launched by the French Government
and implemented by ANR, with the reference ANR‐18‐IdEx‐0001
as part of its program "Emergence".  It was also supported by the Agence Nationale de la Recherche (ANR) Grant ANR-23-CE40-0003 (Project CONVIVIALITY). The author thanks Guillaume Carlier for explaining to him the results of \cite{CFS}, and Patrick Cattiaux, whose work on concentration inequalities inspired this article.

\end{document}